\documentclass[fleqn,letterpaper, 12pt]{article}%
\usepackage{amsfonts,amsmath,latexsym,amssymb,euscript,graphicx}
\usepackage[bookmarks=false,colorlinks=true,linkcolor={blue},pdfstartview={XYZ
null null 1.22}]{hyperref}
\usepackage{graphicx}
\usepackage{color}
\usepackage{amssymb}
\usepackage{amssymb}
\usepackage[T1]{fontenc}
\usepackage{latexsym}
\usepackage{euscript}
\usepackage{amsfonts,amsmath}
\usepackage{verbatim}
\usepackage[english]{babel}
\usepackage{amsmath}
\usepackage{amsfonts}%
\newtheorem{prop}{Proposition}[section]
\newtheorem{cor}[prop]{Corollary}

\newtheorem{lemme}[prop]{Lemma}

\newtheorem{rem}[prop]{Remark}
\newtheorem{remarks}[prop]{Remarks}
\newtheorem{thm}[prop]{Theorem}

\newtheorem{defi}[prop]{Definition}

\renewcommand{\geq}{\geqslant}
\def\leq{\leqslant}

\def\1{{\mathbf{1}}}

\def\1{{\mathbf{1}}}
\def\0.5{{\frac{1}{2}}}

\newcommand{\qed}{\nopagebreak\hspace*{\fill}
{\vrule width6pt height6ptdepth0pt}\par}

\ifx\pdfoutput\relax\let\pdfoutput=\undefined\fi
\newcount\msipdfoutput
\ifx\pdfoutput\undefined\else
\ifcase\pdfoutput\else
\ifx\paperwidth\undefined\else
\ifdim\paperheight=0pt\relax\else\pdfpageheight\paperheight\fi
\ifdim\paperwidth=0pt\relax\else\pdfpagewidth\paperwidth\fi
\fi\fi\fi
\begin{document}

\begin{center}
{\LARGE \textbf{Comparison inequalities on Wiener space}}\newline~\newline
Ivan Nourdin\footnote{Institut Elie Cartan de Lorraine, Universit\'{e} de
Lorraine, BP 70239, 54506 Vandoeuvre-l\`{e}s-Nancy Cedex, France. Email:
\texttt{ivan.nourdin@univ-lorraine.fr}; IN was supported in part by the
(french) ANR grant `Malliavin, Stein and Stochastic Equations with Irregular
Coefficients' [ANR-10-BLAN-0121]}, Giovanni Peccati\footnote{Facult\'{e} des
Sciences, de la Technologie et de la Communication; UR en Math\'{e}matiques.
6, rue Richard Coudenhove-Kalergi, L-1359 Luxembourg, Email:
giovanni.peccati@gmail.com} and Frederi G. Viens\footnote{Dept. Statistics and
Dept. Mathematics, Purdue University, 150 N. University St., West Lafayette,
IN 47907-2067, USA, \texttt{viens@purdue.edu.} FV's research was partially
supported by NSF grant DMS 0907321.}\newline~\newline
\end{center}

{\small \noindent\textbf{Abstract:} We define a covariance-type operator on
Wiener space: for }${\small F}$ {\small and }${\small G}$ {\small two random
variables in the Gross-Sobolev space }${\small D}^{1,2}$ {\small of random
variables with a square-integrable Malliavin derivative, we let }%
${\small \Gamma}_{F,G}{\small :=}\left\langle {\small DF,-DL}^{-1}%
{\small G}\right\rangle $ {\small where }${\small D}$ {\small is the Malliavin
derivative operator and }${\small L}^{-1}$ {\small is the pseudo-inverse of
the generator of the Ornstein-Uhlenbeck semigroup. We use }${\small \Gamma}$
{\small to extend the notion of covariance and canonical metric for vectors
and random fields on Wiener space, and prove corresponding non-Gaussian
comparison inequalities on Wiener space, which extend the Sudakov-Fernique
result on comparison of expected suprema of Gaussian fields, and the Slepian
inequality for functionals of Gaussian vectors. These results are proved using
a so-called smart-path method on Wiener space, and are illustrated via various
examples. We also illustrate the use of the same method by proving a
Sherrington-Kirkpatrick universality result for spin systems in correlated and
non-stationary non-Gaussian random media.

\vspace*{0.15in}

\vspace*{0.15in}


\section{Introduction}

The canonical metric of a centered field $G$ on an index set $T$ is the square
root of the quantity $\delta_{G}^{2}\left(  s,t\right)  =\mathbf{E}\left[
\left(  G_{t}-G_{s}\right)  ^{2}\right]  $, $s,t\in T$. When $G$ is Gaussian,
this $\delta^{2}$ characterizes much of $G\,$'s distribution, and is useful in
various contexts for estimating $G$'s behavior, from its modulus of
continuity, to its expected supremum; see \cite{adler} for an introduction.
The canonical metric, together with the variances of $G$, are of course
equivalent to the covariance function $Q_{G}\left(  s,t\right)  =\mathbf{E}%
\left[  G_{t}G_{s}\right]  $, which defines $G$'s law when $G$ is Gaussian. In
this article, we concentrate on comparison results for expectations of suprema
and other types of functionals, beyond the Gaussian context, by using an
extension of the concepts of covariance and canonical metric on Wiener space.
We introduce these concepts now. For the details of analysis on Wiener space
needed for the next definitions, including the space $\mathbb{D}^{1,2}$ and
the operators $D$ and $L^{-1}$, see Chapter 1 in \cite{Nbook} or Chapter 2 in
\cite{NPbook}. The notion of a `separable random field' is formally defined e.g. in \cite[p. 8]{adlertaylor}.

\begin{defi}
\label{defs}Consider an isonormal Gaussian process $W$ defined on the probability space $\left(  \Omega,\mathcal{F},\mathbf{P}\right) $, and associated with the real separable
Hilbert space $\mathfrak{H}$: recall that this means that $W = \left\{  W\left(  h\right)
:h\in\mathfrak{H}\right\}  $ is a centered Gaussian family such that
$\mathbf{E}\left[  W\left(  h\right)  W\left(  k\right)  \right]
=\left\langle h,k\right\rangle _{\mathfrak{H}}$. Let $\mathbb{D}^{1,2}$ be the
Gross-Sobolev space of random variables $F$ with a square-integrable Malliavin
derivative, i.e. such that $DF\in L^{2}\left(  \Omega\times\mathfrak{H}\right)
$. We denote the generator of the associated Ornstein-Uhlenbeck operator
by $L$.
\noindent For a pair of random variables $F,G\in\mathbb{D}^{1,2}$, we define a
covariance-type operator by%
\begin{equation}
\Gamma_{F,G}:=\langle DF,-DL^{-1}G\rangle_{\mathfrak{H}}. \label{GammaWiener}%
\end{equation}
\noindent Let $F=\{F_{t}\}_{t\in T}$ be a separable random field on an index set $T$,
such that $F_{t}\in\mathbb{D}^{1,2}$ for each $t\in T$. The analogue for the
operator $\Gamma$ of the covariance of $F$ is denoted by%
\begin{equation}
\Gamma_{F}\left(  s,t\right)  :=\Gamma_{F_{s},F_{t}}=\langle D(F_{t}%
),-DL^{-1}(F_{s})\rangle_{\mathfrak{H}}. \label{CovWiener}%
\end{equation}
The analogue for $\Gamma$ of the canonical metric $\delta^{2}$ of $F$ is
denoted by%
\begin{equation}
\Delta_{F}\left(  s,t\right)  :=\langle D(F_{t}-F_{s}),-DL^{-1}(F_{t}%
-F_{s})\rangle_{\mathfrak{H}}. \label{DeltaWiener}%
\end{equation}
\end{defi}

\begin{rem}
\begin{itemize}
\item[(i)] When $F=\{F_{t}\}_{t\in T}$ is in the first Wiener chaos, and hence is
a centered Gaussian field, $\Gamma_{F}$ coincides with its covariance function
$Q_{F}$.
\item[(ii)] In general, the random variable $\Delta_{F}\left(  s,t\right)$ is not positive. However, according e.g. to \cite[Proposition 3.9]{nourdin-peccati}, one has that ${\bf E}[\Delta_F(s,t) | F_t-F_s]\geq 0$, a.s.-${\bf P}$. 
\end{itemize}
\end{rem}

The extension of the concept of covariance function given above in
(\ref{GammaWiener}) appeared in \cite{AMV} and in \cite{NPR}, respectively to
aid in the study of densities of random vectors and of multivariate normal
approximations, both on Wiener space. Comparison results on Wiener space have,
in the past, focused on concentration or Poincar\'{e} inequalities: see
\cite{U}. Recently, the scalar analogue of the covariance operator above, i.e.
$\Gamma_{F,F}$, was exploited to derive sharp tail comparisons on Wiener
space, in \cite{nourdin-viens} and \cite{V}.

The two main types of comparison results we will investigate herein are those
of Sudakov-Fernique type and those of Slepian type. See \cite{adler, adlertaylor} for
details of the classical proofs.

In the basic \emph{Sudakov-Fernique} inequality, one considers two centered
separable Gaussian fields $F$ and $G$ on $T$, such that $\delta_{F}^{2}\left(
s,t\right)  \geq\delta_{G}^{2}\left(  s,t\right)  $ for all $s,t\in T$; then
$\mathbf{E}\left[  \sup_{T}F\right]  \geq\mathbf{E}\left[  \sup_{T}G\right]
$. Here $T$ can be any index set, as long as the laws of $F$ and $G$ can be
determined by considering only countably many elements of $T$; this works for
instance if $T$ is a subset of Euclidean space and $F$ and $G$ are a.s.
continuous. To try to extend this result to non-Gaussian fields with no
additional machinery, for illustrative purposes, the following setup provides
an easy example.

\begin{prop}
\label{sad-prop} Let $F\ $and $G\ $be two separable fields on $T$, with $G$
and $F-G$ independent, and $E[F_{t}]=E[G_{t}]$ for every $t\in T$. Then
$\mathbf{E}\left[  \sup_{T}F\right]  \geq\mathbf{E}\left[  \sup_{T}G\right]  $.
\end{prop}

The proof of this proposition is elementary. Let $H=F-G$. Note that for any
$t_{0}\in T$, $\mathbf{E}\left[  H\left(  t_{0}\right)  \right]  =0$. We may
write $\mathbf{P=P}_{H}\times\mathbf{P}_{F}$ with obvious notation. Thus
\[
\mathbf{E}\left[  \sup_{T}F\right]  =\mathbf{E}\left[  \sup_{T}\left(
H+G\right)  \right]  =\mathbf{E}_{G}\left[  \mathbf{E}_{H}\left[  \sup
_{T}\left(  H+G\right)  \right]  \right]
\]
where under $\mathbf{P}_{H}$, $G$ is deterministic. Thus
\[
\mathbf{E}\left[  \sup_{T}F\right]  \geq\mathbf{E}_{G}\left[  \mathbf{E}%
_{H}\left[  H\left(  t_{0}\right)  +\sup_{T}G\right]  \right]  =\mathbf{E}%
_{G}\left[  \mathbf{E}_{H}\left[  H\left(  t_{0}\right)  \right]  +\sup
_{T}G\right]  =\mathbf{E}_{G}\left[  \sup_{T}G\right] .
\]

What makes this proposition so easy to establish is the very strong joint
distributional assumption on $\left(  F,G\right)  $, even though we do not
make any marginal distributional assumptions about $F$ and $G$. Also note that
in the Gaussian case, the covariance assumption on $\left(  F,G\right)  $
implies that $\delta_{F}^{2}\left(  s,t\right)  \geq\delta_{G}^{2}\left(
s,t\right)  $, and is in fact a much stronger assumption than simply comparing
these canonical metrics, so that the classical Sudakov-Fernique inequality
applies handily.

Let us now discuss the \emph{Slepian inequality} similarly. In the basic
inequality, consider two centered Gaussian vectors $F$ and $G$ in
$\mathbb{R}^{d}$, with covariance matrices $\left(  B_{ij}\right)  $ and
$\left(  C_{ij}\right)  $. Let $f\in C^{2}\left(  \mathbb{R}^{d}\right)  $
with bounded partial derivatives up to order 2. Assume that for all
$x\in\mathbb{R}^{d}$,%
\[
\sum_{i,j=1}^{d}\left(  B_{ij}-C_{ij}\right)  \frac{\partial^{2}f}{\partial
x_{i}\partial x_{j}}(x)\geq0.
\]
Then $\mathbf{E}\left[  f\left(  F\right)  \right]  \geq\mathbf{E}\left[
f\left(  G\right)  \right]  $. To obtain such a result for non-Gaussian
vectors, one may again try to impose strong joint-distributional conditions to
avoid marginal conditions. The following example is a good illustration. With
$F$ and $G$ two random vectors in $\mathbb{R}^{d}$ and $f$ convex on
$\mathbb{R}^{d}$, assume that $\mathbf{E}[F]=\mathbf{E}[G]$, $\mathbf{E}%
|f(F)|<\infty$, $\mathbf{E}|f(G)|<\infty$, and $G$ and $F-G$ are independent.
By convexity for any $c\in\mathbb{R}^{d}$ we have that
\[
f(F-G+c)\geq f(c)+\langle\nabla f(c),F-G\rangle_{\mathbb{R}^{d}}.
\]
Hence $\mathbf{E}[f(F-G+c)]\geq f(c)$. By choosing $c=G$ and then taking
expectations, we get $\mathbf{E}[f(F)]\geq\mathbf{E}[f(G)]$, i.e. the Slepian
inequality conclusion holds. In other word we have the following.

\begin{prop}
\label{gordon-prop}Let $F\ $and $G\ $be two random vectors in $\mathbb{R}^{d}%
$, with $G$ and $F-G$ independent. Let $f:\mathbb{R}^{d}\to\mathbb{R}$ be a
convex function. Assume $\mathbf{E}[F]=\mathbf{E}[G]$, $\mathbf{E}%
|f(F)|<\infty$, $\mathbf{E}|f(G)|<\infty$. Then $\mathbf{E}\left[
f(F)\right]  \geq\mathbf{E}\left[  f(G)\right]  $.
\end{prop}

To avoid very strong joint law assumptions on $\left(  F,G\right)  $ such as
those used in the two elementary propositions above, this paper concentrates
instead on exploiting some mild assumptions on the marginals of $F$ and $G$,
particularly imposing Malliavin differentiability as in Definition \ref{defs}.
We will see in particular that to obtain a Sudakov-Fernique inequality for
highly non-Gaussian fields, one can use $\Delta$ instead of $\delta^{2}$, and
to get a Slepian inequality in the same setting, one can use $\Gamma
_{F_{i},F_{j}}$ and $\Gamma_{G_{i},G_{j}}$ instead of $B_{i,j}$ and $C_{i,j}$
respectively. The proofs we use are based on the technique of interpolation,
and on the following integration-by-parts theorem on Wiener space, which was
first introduced in \cite{nourdin-peccati} (also see Theorem 2.9.1 in
\cite{NPbook}): for any centered $F,G\in\mathbb{D}^{1,2}$, $\mathbf{E}\left[
FG\right]  =\mathbf{E}\left[  \Gamma_{F,G}\right]  .$ This formula is
particularly useful when combined with the chain rule of the Malliavin
calculus, to yield that for any $\Phi:\mathbb{R\to R}$ such that
$\mathbf{E}\left[  \Phi^{\prime}\left(  F\right)  ^{2}\right]  <\infty$,
\begin{equation}
\mathbf{E}\left[  \Phi\left(  F\right)  G\right]  =\mathbf{E}\left[
\Phi^{\prime}\left(  F\right)  \Gamma_{F,G}\right]  . \label{parts}%
\end{equation}

The remainder of this paper is structured as follows. In Section \ref{SUP}, we
prove a new Sudakov-Fernique inequality for comparing suprema of random fields
on Wiener space, and show how this may be applied to the supremum of the
solution of a stochastic differential equation with non-linear drift, driven
by fractional Brownian motion. In Section \ref{GS}, we prove a Slepian-type
inequality for comparing non-linear functionals of random vectors on Wiener
space, and apply it to a comparison result for perturbations of Gaussian
vectors, and to a concentration inequality. Finally in Section \ref{SK}, we
show how to extend the universality class of the Sherrington-Kirkpatrick spin
system, to some random media on Wiener space with dependence and
non-stationarity. All our main theorems' proofs are based on the extension to
Wiener space of the so-called smart-path method using the objects identified
in Definition \ref{defs}.

\section{A result of Sudakov-Fernique type\label{SUP}}

The proof of the following result is based on an extension of classical computations based on a `smart path method' that are available in the Gaussian setting. The reader is referred to \cite[p. 61]{adlertaylor} for a similar proof (originally due to S. Chatterjee, see also \cite{Chatterjee}) in the simpler Gaussian setting.

\begin{thm}
\label{fs-thm} Let $F=\{F_{t}\}_{t\in T}$ and $G=\{G_{t}\}_{t\in T}$ be
separable centered random fields on an index set $T$, such that $F_{t}%
,G_{t}\in\mathbb{D}^{1,2}$ for every $t\in T$. Their canonical metrics on
Wiener space, $\Delta_{F}$ and $\Delta_{G}$, are defined according to
(\ref{DeltaWiener}). Assume that $\mathbf{E}\left[  {\sup_{T}}F\right]
<\infty$ and $\mathbf{E}\left[  {\sup_{T}}G\right]  <\infty$. Assume that
almost surely for all $s,t\in T$,
\begin{equation}
\Delta_{F}\left(  s,t\right)  \leqslant\Delta_{G}\left(  s,t\right)  .
\label{assu1}%
\end{equation}
Assume furthermore that almost surely for all $s,t\in T$,%
\begin{equation}
\Gamma_{F_{s},G_{t}}=0. \label{assu2}%
\end{equation}
Then%
\[
\mathbf{E}\left[  \displaystyle{\sup_{t\in T}}\,F_{t}\right]  \leqslant
\mathbf{E}\left[  \displaystyle{\sup_{t\in T}}\,G_{t}\right]  .
\]

\end{thm}

\begin{rem}
\label{SFV}\textrm{If }$\left(  F,G\right)  $ \textrm{is jointly Gaussian, one
can assume that both processes belong to the first Wiener chaos, and then}
\[
\langle D(F_{t}-F_{s}),-DL^{-1}(F_{t}-F_{s})\rangle_{\mathfrak{H}}%
=E[(F_{t}-F_{s})^{2}],
\]
\textrm{and similarly for }$G$\textrm{. The orthogonality condition
(\ref{assu2}) is then equivalent to independence. As such, Theorem
\ref{fs-thm} extends the classical Sudakov-Fernique inequality, as stated e.g.
in Vitale \cite[Theorem 1]{Vitale} in the case $|$}$T|<\infty$\textrm{.}
\end{rem}

\begin{cor}
\label{fs-cor1}When $G$ belongs to the first Wiener chaos (in particular, $G$
is Gaussian), then $\Delta_{G}\left(  s,t\right)  =\delta_{G}^{2}\left(
s,t\right)  $ is $G$'s (non-random) canonical metric, and the conclusion of
Theorem \ref{fs-thm} continues to hold without Assumption (\ref{assu2}).
\end{cor}

\textit{Proof}. Since $\Delta_{G}$ is non-random, the Gaussian process $G$ in
this corollary can be defined on any probability space, and thus we can assume
that $G$ is independent of $F$, and therefore that Assumption (\ref{assu2})
holds.\nopagebreak\hspace*{\fill} {\vrule width6pt height6ptdepth0pt}\bigskip

 \textit{Proof of Theorem \ref{fs-thm}}.

 $\emph{Step}$\emph{ 1: Approximation}. For each $n>0$, let $T_{n}$ be a finite
subset of $T$ such that $T_{n}\subset T_{n+1}$ and $T_{n}$ increases to a
countable subset of $T$ on which the laws of $F$ and $G$ are determined (for
instance, if $T=\mathbb{R}_{+}$ and $F$ and $G$ are continuous, we may choose
for $T_{n}$ the set of dyadics of order $n$). By separability, as
$n\rightarrow\infty$,
\[
\sup_{t\in T_{n}}F_{t}\overset{\mathrm{a.s.}}{\rightarrow}\sup_{t\in T}%
F_{t}\quad\mbox{and}\quad\sup_{t\in T_{n}}G_{t}\overset{\mathrm{a.s.}%
}{\rightarrow}\sup_{t\in T}G_{t}%
\]
and, since the convergence is monotone, we also have that as $n\rightarrow
\infty$,%
\[
\mathbf{E}\left[  \sup_{t\in T_{n}}F_{t}\right]  \rightarrow\mathbf{E}\left[
\sup_{t\in T}F_{t}\right]  \quad\mbox{and}\quad\mathbf{E}\left[  \sup_{t\in
T_{n}}F_{t}\right]  \rightarrow\mathbf{E}\left[  \sup_{t\in T}F_{t}\right]  .
\]
Therefore, we assume without loss of generality in the remainder of the proof
that $T=\{1,2,\ldots,d\}$ is finite.\bigskip

 \emph{Step 2: calculation.} Fix $\beta>0$, and consider, for any $t\in
\lbrack0,1]$,
\[
\varphi(t)=\frac{1}{\beta}\,\mathbf{E}\left[  \log\left(  \sum_{i=1}%
^{d}e^{\beta(\sqrt{1-t}G_{i}+\sqrt{t}F_{i})}\right)  \right]  .
\]
Let us differentiate $\varphi$ with respect to $t\in(0,1)$. We get%
\begin{equation}
\varphi^{\prime}(t)=\frac{1}{2}\sum_{i=1}^{d}\mathbf{E}\left[  \left(
\frac{1}{\sqrt{t}}F_{i}-\frac{1}{\sqrt{1-t}}G_{i}\right)  h_{t,\beta
,i}(F,G)\right]  , \label{phiprime}%
\end{equation}
where, for $x,y\in\mathbb{R}^{d}$, $i=1,\ldots,d$, $t\in(0,1)$ and $\beta>0$,
we set
\[
h_{t,\beta,i}(x,y)=\frac{e^{\beta(\sqrt{1-t}y_{i}+\sqrt{t}x_{i})}}{\sum
_{j=1}^{d}e^{\beta(\sqrt{1-t}y_{j}+\sqrt{t}x_{j})}}.
\]
Using the integration-by-parts formula (\ref{parts}) in (\ref{phiprime})
yields%
\begin{align}
&  \varphi^{\prime}(t)\nonumber\\
&  =\frac{1}{2}\sum_{i,j=1}^{d}\left(  \frac{1}{\sqrt{t}}\mathbf{E}\left[
\frac{\partial h_{t,\beta,i}}{\partial x_{j}}(F,G)\Gamma_{F_{j},F_{i}}\right]
-\frac{1}{\sqrt{1-t}}\mathbf{E}\left[  \frac{\partial h_{t,\beta,i}}{\partial
y_{j}}(F,G)\Gamma_{G_{j},G_{i}}\right]  \right) \nonumber\\
&  +\frac{1}{2}\sum_{i,j=1}^{d}\left(  \frac{1}{\sqrt{t}}\mathbf{E}\left[
\frac{\partial h_{t,\beta,i}}{\partial x_{j}}(F,G)\Gamma_{G_{j},F_{i}}\right]
-\frac{1}{\sqrt{1-t}}\mathbf{E}\left[  \frac{\partial h_{t,\beta,i}}{\partial
y_{j}}(F,G)\Gamma_{F_{j},G_{i}}\right]  \right)  . \label{phiprime2}%
\end{align}
The orthogonality assumption (\ref{assu2}) implies that all the terms in the
last line of (\ref{phiprime2}) are zero. For $i\neq j$, we have
\begin{align*}
\frac{\partial h_{t,\beta,i}}{\partial x_{i}}(x,y)  &  =\beta\sqrt
{t}\big(h_{t,\beta,i}(x,y)-h_{t,\beta,i}(x,y)^{2}\big)\\
\frac{\partial h_{t,\beta,i}}{\partial x_{j}}(x,y)  &  =-\beta\sqrt
{t}\,h_{t,\beta,i}(x,y)h_{t,\beta,j}(x,y)\\
\frac{\partial h_{t,\beta,i}}{\partial y_{i}}(x,y)  &  =\beta\sqrt
{1-t}\big(h_{t,\beta,i}(x,y)-h_{t,\beta,i}(x,y)^{2}\big)\\
\frac{\partial h_{t,\beta,i}}{\partial y_{j}}(x,y)  &  =-\beta\sqrt
{1-t}\,h_{t,\beta,i}(x,y)h_{t,\beta,j}(x,y).
\end{align*}
Therefore
\begin{align*}
\varphi^{\prime}(t)  &  =\frac{\beta}{2}\sum_{i}\mathbf{E}\bigg[h_{t,\beta
,i}(F,G)(1-h_{t,\beta,i}(F,G))\big(\Gamma_{F_{i},F_{i}}-\Gamma_{G_{i},G_{i}%
}\big)\bigg]\\
&  -\frac{\beta}{2}\sum_{i\neq j}\mathbf{E}\bigg[h_{t,\beta,i}(F,G)h_{t,\beta
,j}(F,G)\big(\Gamma_{F_{i},F_{j}}-\Gamma_{G_{i},G_{j}}\big)\bigg]\\
&  =\frac{\beta}{2}\sum_{i}\mathbf{E}\bigg[h_{t,\beta,i}(F,G)\big(\Gamma
_{F_{i},F_{i}}-\Gamma_{G_{i},G_{i}}\big)\bigg]\\
&  -\frac{\beta}{2}\sum_{i,j}\mathbf{E}\bigg[h_{t,\beta,i}(F,G)h_{t,\beta
,j}(F,G)\big(\Gamma_{F_{i},F_{j}}-\Gamma_{G_{i},G_{j}}\big)\bigg].
\end{align*}
But $\sum_{i=1}^{d}h_{t,\beta,i}(F,G)=1$, hence $\varphi^{\prime}(t)$ is given
by%
\[
\frac{\beta}{4}\sum_{i,j=1}^{d}\mathbf{E}\bigg[h_{t,\beta,i}(F,G)h_{t,\beta
,j}(F,G)\big(\Delta_{F}\left(  i,j\right)  -\Delta_{G}\left(  i,j\right)
\big)\bigg].
\]
\bigskip

 \emph{Step 3: estimation and conclusion. }We observe that $h_{t,\beta
,i}(F,G)>0$ for all $i$. Moreover, by assumption (\ref{assu1}) we get
$\varphi^{\prime}(t)\leqslant0$ for all $t$, implying in turn that
$\varphi(0)\geq\varphi(1)$, that is
\[
\frac{1}{\beta}\,\mathbf{E}\left[  \log\left(  \sum_{i=1}^{d}e^{\beta F_{i}%
}\right)  \right]  \leqslant\frac{1}{\beta}\,\mathbf{E}\left[  \log\left(
\sum_{i=1}^{d}e^{\beta G_{i}}\right)  \right]
\]
for any $\beta>0$. But
\[
\displaystyle{\max_{1\leqslant i\leqslant d}}F_{i}=\frac{1}{\beta}%
\,\log\left(  e^{\beta\times\displaystyle{\max_{1\leqslant i\leqslant d}}%
F_{i}}\right)  \leqslant\frac{1}{\beta}\log\left(  \sum_{i=1}^{d}e^{\beta
F_{i}}\right)  \leqslant\frac{\log d}{\beta}+\displaystyle{\max_{1\leqslant
i\leqslant d}}F_{i},
\]
and the same with $G$ instead of $F$. Therefore
\[
\mathbf{E}\left[  \displaystyle{\max_{1\leqslant i\leqslant d}}F_{i}\right]
\leqslant\mathbf{E}\left[  \frac{1}{\beta}\log\left(  \sum_{i=1}^{d}e^{\beta
F_{i}}\right)  \right]  \leqslant\mathbf{E}\left[  \frac{1}{\beta}\log\left(
\sum_{i=1}^{d}e^{\beta G_{i}}\right)  \right]  \leqslant\frac{\log d}{\beta
}+\mathbf{E}\left[  \displaystyle{\max_{1\leqslant i\leqslant d}}G_{i}\right]
,
\]
and the desired conclusion follows by letting $\beta$ goes to infinity.
\nopagebreak\hspace*{\fill} {\vrule width6pt height6ptdepth0pt}\bigskip

We now give an example of application of Theorem \ref{fs-thm}, to a problem of
current interest in stochastic analysis.

\subsection{Example: supremum of an SDE driven by fBm}

Let $B^{H}$ be a fractional Brownian motion with Hurst index $H>1/2$, let
$b:\mathbb{R}\rightarrow\mathbb{R}$ be increasing and Lipschitz (in
particular, $b^{\prime}\geq0$ almost everywhere), and let $x_{0}\in\mathbb{R}%
$. We consider the process $F=(F_{t})_{t\in\lbrack0,T]}$ defined as the unique
solution to
\begin{equation}
F_{t}=x_{0}+B_{t}^{H}+\int_{0}^{t}b(F_{s})ds. \label{edsfrac}%
\end{equation}
(For more details about this equation, we refer the reader to
\cite{nualart-ouknine}.) It is well-known (see e.g. \cite{nualart-saussereau}
or \cite{nourdin-simon}) that, for any $t\in(0,T]$, we have that $F_{t}%
\in\mathbb{D}^{1,2}$ with
\begin{equation}
D_{u}F_{t}=\mathbf{1}_{[0,t]}(u)\,\exp\left(  \int_{u}^{t}b^{\prime}%
(F_{w})dw\right)  . \label{duft}%
\end{equation}
Fix $t>s\geq0$. By combining (\ref{duft}) with a calculation technique
described e.g. in \cite[Proposition 3.7]{nourdin-viens} based on the so-called
Mehler formula, we get
\begin{align}
&  \Delta_{F}\left(  s,t\right) \label{fbm1}\\
&  =H(2H-1)\widehat{E}\bigg\{\int_{0}^{\infty}e^{-z}\bigg[\int_{[0,s]^{2}%
}\left(  e^{\int_{u}^{t}b^{\prime}(F_{w})dw}-e^{\int_{u}^{s}b^{\prime}%
(F_{w})dw}\right) \nonumber\\
&  \hskip5cm\times\left(  e^{\int_{v}^{t}b^{\prime}(F_{w}^{(z)})dw}%
-e^{\int_{v}^{s}b^{\prime}(F_{w}^{(z)})dw}\right)  |u-v|^{2H-2}dudv\nonumber\\
&  +\int_{[0,s]\times\lbrack s,t]}\left(  e^{\int_{u}^{t}b^{\prime}(F_{w}%
)dw}-e^{\int_{u}^{s}b^{\prime}(F_{w})dw}\right)  \,e^{\int_{v}^{t}b^{\prime
}(F_{w}^{(z)})dw}\,|u-v|^{2H-2}dudv\nonumber\\
&  +\int_{[s,t]\times\lbrack0,s]}e^{\int_{u}^{t}b^{\prime}(F_{w})dw}\left(
e^{\int_{v}^{t}b^{\prime}(F_{w}^{(z)})dw}-e^{\int_{v}^{s}b^{\prime}%
(F_{w}^{(z)})dw}\right)  |u-v|^{2H-2}dudv\nonumber\\
&  +\int_{[s,t]^{2}}e^{\int_{u}^{t}b^{\prime}(F_{w})dw+\int_{v}^{t}b^{\prime
}(F_{w}^{(z)})dw}|u-v|^{2H-2}dudv\bigg]dz\bigg\}.\nonumber
\end{align}
Here, $F^{(z)}$ means the solution to (\ref{edsfrac}), but when $B^{H}$ is
replaced by the new fractional Brownian motion $e^{-z}B^{H}+\sqrt{1-e^{-2z}%
}\widehat{B}^{H}$, for $\widehat{B}^{H}$ an independent copy of $B^{H}$, and
$\widehat{E}$ is the mathematical expectation with respect to $\widehat{B}%
^{H}$ only. Because $b^{\prime}\geq0$, we see that
\begin{align*}
\exp\left\{  \int_{u}^{t}b^{\prime}(F_{w})dw\right\}  -\exp\left\{  \int
_{u}^{s}b^{\prime}(F_{w})dw\right\}   &  \geq0\quad
\mbox{for any $0\leq u\leq s<t$,}\\
\exp\left\{  \int_{v}^{t}b^{\prime}(F_{w}^{(z)})dw\right\}  -\exp\left\{
\int_{v}^{s}b^{\prime}(F_{w}^{(z)})dw\right\}   &  \geq0\quad
\mbox{for any $0\leq v\leq s<t$,}\\
\exp\left\{  \int_{u}^{t}b^{\prime}(F_{w})dw+\int_{v}^{t}b^{\prime}%
(F_{w}^{(z)})dw\right\}   &  \geq1\quad\mbox{for any $s\leq u,v\leq t$}.
\end{align*}
In particular, $\Delta_{F}\left(  s,t\right)  \geq H(2H-1)\int_{[s,t]^{2}%
}|u-v|^{2H-2}dudv=|t-s|^{2H}$. We recognize $|t-s|^{2H}$ as the squared
canonical metric of fractional Brownian motion, and we deduce from Theorem
\ref{fs-thm} (observe that it is not a loss of generality to have assumed that
$s<t$) that
\[
E\left[  \displaystyle{\max_{t\in\lbrack0,T]}}\big(F_{t}-E[F_{t}]\big)\right]
\geq E\left[  \displaystyle{\max_{t\in\lbrack0,T]}}B_{t}^{H}\right]  .
\]
Also note that by the same calculation as above, the inequality in the
conclusion is reversed if $b$ is decreasing.

\section{A result of Slepian type\label{GS}}

In Section \ref{SUP}, we investigated the ability to compare suprema of random
vectors and fields based on covariances and the Wiener-space extensions of the
concept of covariance in Definition \ref{defs}. In this section, we show that
these extensions also apply to functionals beyond the supremum, under
appropriate convexity assumptions.

\begin{thm}
\label{gs-thm} Let $F,G$ be two centered rv's in $\mathbb{D}^{1,2}\left(
\mathbb{R}^{d}\right)  $, in other words, assume that for every $i=1,2,\cdots
,d$, $F_{i}\in\mathbb{D}^{1,2}$ and $G_{i}\in\mathbb{D}^{1,2}$ and
$\mathbf{E}[F_{i}]=\mathbf{E}[G_{i}]=0$. Let also $f:\mathbb{R}^{d}%
\rightarrow\mathbb{R}$ be a $C^{2}$-function. We define the $d\times d$ random
\textquotedblleft covariance\textquotedblright-type matrix $$\Gamma
^{F}=\left\{  \Gamma_{ij}^{F}:=\Gamma_{F_{i},F_{j}}:i,j=1,\cdots,d\right\}  $$
for $F$, according to (\ref{GammaWiener}), and similarly for $\Gamma^{G}$. We
assume that $\mathbf{E}\left[  \left\vert \frac{\partial^{2}f}{\partial
x_{i}\partial x_{j}}(\sqrt{1-t}G+\sqrt{t}F)\right\vert \right]  $ is finite
for every $i,j =1,\cdots,d$ and $t\in\lbrack0,1]$, that $\Gamma_{F_{i},G_{j}%
}=0$ for any $i,j$ and that for all $x\in\mathbb{R}^{d}$, almost surely,
\begin{equation}
\sum_{i,j=1}^{d}\left(  \Gamma_{ij}^{F}-\Gamma_{ij}^{G}\right)  \frac
{\partial^{2}f}{\partial x_{i}\partial x_{j}}(x)\geq0. \label{gs-condition}%
\end{equation}
Then $\mathbf{E}[f(F)]\geq\mathbf{E}[f(G)]$.
\end{thm}

\begin{rem}
\textrm{If }$F$\textrm{ and }$G$\textrm{ are Gaussian, then }$\Gamma^{F}%
$\textrm{ and }$\Gamma^{G}$ \textrm{are the covariance matrices of }%
$F$\textrm{ and }$G$\textrm{ a.s.},\textrm{ and we recover the classical
Slepian inequality, see e.g. \cite{slepian}, or the paragraph in the
Introduction preceding Proposition \ref{gordon-prop}. }
\end{rem}

\begin{cor}
If $F$ is Gaussian (but not necessarily $G$), then the conclusion of Theorem
\ref{gs-thm} holds without any information on the joint law of $\left(
F,G\right)  $, except for assuming that if $F$ and $G$ are independent, then
$\mathbf{E}\left[  \left\vert \frac{\partial^{2}f}{\partial x_{i}\partial
x_{j}}(\sqrt{1-t}G+\sqrt{t}F)\right\vert \right]  $ is finite for every $i,j
=1,\cdots,d$ and $t\in\lbrack0,1]$.
\end{cor}

\textit{Proof of Theorem \ref{gs-thm}}. For $t\in\lbrack0,1]$, set
\[
\varphi(t)=E[f(\sqrt{1-t}G+\sqrt{t}F)].
\]
We have
\[
\varphi^{\prime}(t)=\frac{1}{2}\sum_{i=1}^{d}\left(  \frac{1}{\sqrt{t}%
}E\left[  \frac{\partial f}{\partial x_{i}}(\sqrt{1-t}G+\sqrt{t}%
F)F_{i}\right]  -\frac{1}{\sqrt{1-t}}E\left[  \frac{\partial f}{\partial
x_{i}}(\sqrt{1-t}G+\sqrt{t}F)G_{i}\right]  \right)  .
\]
By using the integrating-by-parts formula (\ref{parts}), we get the following
extension of a classical identity due to Piterbarg \cite{piterbarg}:
\begin{align*}
\varphi^{\prime}(t)  &  =\frac{1}{2}\sum_{i,j=1}^{d}E\left[  \frac
{\partial^{2}f}{\partial x_{i}\partial x_{j}}(\sqrt{1-t}G+\sqrt{t}%
F)\big(\langle DF_{j},-DL^{-1}F_{i}\rangle_{\mathfrak{H}}-\langle
DG_{j},-DL^{-1}G_{i}\rangle_{\mathfrak{H}}\big)\right] \\
&  =\frac{1}{2}\sum_{i,j=1}^{d}E\left[  \frac{\partial^{2}f}{\partial
x_{i}\partial x_{j}}(\sqrt{1-t}G+\sqrt{t}F)\left(  \Gamma_{ij}^{F}-\Gamma
_{ij}^{G}\right)  \right]  .
\end{align*}
As a consequence, $\varphi^{\prime}(t)\geq0$ (resp. $\leqslant$), implying in
turn that $\varphi(1)\geq\varphi(0)$ (resp. $\leqslant$), which is the desired
conclusion. \nopagebreak\hspace*{\fill} {\vrule width6pt height6ptdepth0pt}%
\vspace*{0.15in}

\textit{Proof of the corollary.} When $F$ is Gaussian, $\Gamma^{F}$ is almost
surely deterministic, and we may thus assume that $F$ and $G$ are defined on
the same probability space and are independent. The finiteness of
$\mathbf{E}\left[  \left\vert \frac{\partial^{2}f}{\partial x_{i}\partial
x_{j}}(\sqrt{1-t}G+\sqrt{t}F)\right\vert \right]  $ can then be assumed to
hold, and the theorem applies. \nopagebreak\hspace*{\fill} {\vrule width6pt
height6ptdepth0pt}

\subsection{Example: perturbation of a Gaussian vector}

Here we present an example of how to perturb an arbitrary Gaussian vector
$G\in\mathbb{R}^{d}$ using a functional on Wiener space to guarantee that for
any function $f$ with non-negative (resp. non-positive) second derivatives,
$f\left(  G\right)  $ sees its expectation increase (resp. decrease) with the
perturbation. It is sufficient for the perturbation to be based on variables
that are positively \textquotedblleft correlated\textquotedblright\ to $G$, in
a sense defined using the covariance operator $\Gamma$ of Definition
\ref{defs}. Let $C$ be the covariance matrix of $G$.

We may assume that for every $i=1,\ldots,d$, $G_{i}=I_{1}\left(  g_{i}\right)
$ where the $g_{i}$'s are such that $\left\langle g_{i},g_{j}\right\rangle
_{\mathfrak{H}}=C_{i,j}$. Fix integers $n_{1},\ldots,n_{d}\geq1$, let
$f_{i,k}$ $i=1,\ldots,d$, $k=1,\ldots,n_{d}$, be a sequence of elements of $H$
such that $\langle f_{i,k},g_{j}\rangle_{\mathfrak{H}}\geq0$ and $\langle
f_{i,k},f_{j,l}\rangle_{\mathfrak{H}}\geq0$ for all $i,j,k,l$, and let
$\Phi_{i}:\mathbb{R}^{n_{i}}\rightarrow\mathbb{R}$, $i=1,\ldots,d$, be a
sequence of $C^{1}$-functions such that $\frac{\partial\Phi_{i}}{\partial
x_{k}}\geq0$ for all $k$ (each $\Phi_{i}$ is increasing w.r.t. every
component). For $i=1,\ldots,d$, we set
\[
F_{i}=G_{i}+\Phi_{i}\big(I_{1}(f_{i,1}),\ldots,I_{1}(f_{i,n_{i}})\big).
\]
Our assumptions are simply saying that all the Gaussian pairs $\left(
G_{j},I_{1}\left(  f_{i,k}\right)  \right)  $ are non-negatively correlated,
as are all the Gaussian pairs $\left(  I_{1}\left(  f_{i,k}\right)
,I_{1}\left(  f_{j,\ell}\right)  \right)  $. For any $i,j=1,\ldots,d$, we
compute
\begin{align*}
DF_{i}  &  =g_{i}+\sum_{k=1}^{n_{i}}\frac{\partial\Phi_{i}}{\partial x_{k}%
}(I_{1}(f_{i,1}),\ldots,I_{1}(f_{i,n_{i}}))f_{i,k}\\
P_{z}DF_{j}  &  =g_{j}+\sum_{l=1}^{n_{j}}\widehat{E}\left[  \frac{\partial
\Phi_{j}}{\partial x_{l}}(I_{1}^{(z)}(f_{j,1}),\ldots,I_{1}^{(z)}(f_{j,n_{j}%
}))\right]  f_{j,l},
\end{align*}
where $I_{1}^{(z)}$ means that the Wiener integral is taken with respect to
$W^{(z)}=e^{-z}W+\sqrt{1-e^{-2z}}\widehat{W}$ instead of $W$, for $\widehat
{W}$ an independent copy of $W$, and where $\widehat{E}$ is the mathematical
expectation with respect to $\widehat{W}$ only. Therefore, using the
Mehler-formula representation of $DL^{-1}$ (see \cite{nourdin-viens}),%
\begin{align*}
&  \Gamma_{i,j}:=\Gamma_{F_{i},F_{j}}=\int_{0}^{\infty}e^{-z}\langle
DF_{i},P_{z}DF_{j}\rangle_{\mathfrak{H}}dz\\
&  =C_{i,j}+\sum_{k=1}^{n_{i}}\frac{\partial\Phi_{i}}{\partial x_{k}}%
(I_{1}(f_{i,1}),\ldots,I_{1}(f_{i,n_{i}}))\langle f_{i,k},g_{j}\rangle
_{\mathfrak{H}}\\
&  +\sum_{l=1}^{n_{j}}\langle f_{j,l},g_{i}\rangle_{\mathfrak{H}}\int
_{0}^{\infty}e^{-z}\widehat{E}\left[  \frac{\partial\Phi_{j}}{\partial x_{l}%
}(I_{1}^{(z)}(f_{j,1}),\ldots,I_{1}^{(z)}(f_{j,n_{j}}))\right]  dz\\
&  +\sum_{k=1}^{n_{i}}\sum_{l=1}^{n_{j}}\langle f_{i,k},f_{j,l}\rangle
_{\mathfrak{H}}\frac{\partial\Phi_{i}}{\partial x_{k}}(I_{1}(f_{i,1}%
),\ldots,I_{1}(f_{i,n_{i}}))\int_{0}^{\infty}e^{-z}\widehat{E}\left[
\frac{\partial\Phi_{j}}{\partial x_{l}}(I_{1}^{(z)}(f_{j,1}),\ldots
,I_{1}^{(z)}(f_{j,n_{j}}))\right]  dz.
\end{align*}
Using the assumptions, we see that $\Gamma_{i,j}\geq\langle g_{i},g_{j}%
\rangle_{\mathfrak{H}}$ for all $i,j=1,\ldots,d$. Hence, for all $C^{2}%
$-function $\Psi:\mathbb{R}^{d}\rightarrow\mathbb{R}$ such that $\frac
{\partial^{2}\Psi}{\partial x_{i}\partial x_{j}}(x)\geq0$ (resp. $\leqslant$),
condition (\ref{gs-condition}) is in order, so that $E[\Psi(F)]\geq
E[\Psi(G)]$ (resp. $\leqslant$) by virtue of Theorem \ref{gs-thm}.

\subsection{Example: a concentration inequality}

Next we encounter an application of Theorem \ref{gs-thm} to compare
distributions of non-Gaussian vectors to Gaussian distributions.

\begin{cor}
\label{Cop-thm}Let $F=(F_{1},\ldots,F_{d})\in\mathbb{R}^{d}$ be such that
$F_{i}\in\mathbb{D}^{1,2}$ and $E[F_{i}]=0$ for every $i$, and define
$\Gamma=\left\{  \Gamma_{ij}:=\Gamma_{F_{i},F_{j}}:i,j=1,\cdots,d\right\}  $,
according to (\ref{GammaWiener}). Let $C\ $be a deterministic non-negative
definite $d\times d$ matrix such that, almost surely, $C-\Gamma$ is
non-negative definite. Then, with $\Vert C\Vert_{op}$ the operator norm of
$C$, for any $x_{1},\ldots,x_{d}\geq0$, we have
\[
P[F_{1}\geq x_{1},\ldots,F_{d}\geq x_{d}]\leqslant\exp\left\{  -\frac
{x_{1}^{2}+\ldots+x_{d}^{2}}{2\Vert C\Vert_{op}}\right\}  .
\]

\end{cor}

\textit{Proof}. For any $\theta\in\mathbb{R}_{+}^{d}$, we can write
\[
P[F_{1}\geq x_{1},\ldots,F_{d}\geq x_{d}]\leqslant P\big[\langle
\theta,F\rangle_{\mathbb{R}^{d}}\geq\langle\theta,x\rangle_{\mathbb{R}^{d}%
}\big]\leqslant e^{-\langle\theta,x\rangle_{\mathbb{R}^{d}}}\,E[e^{\langle
\theta,F\rangle_{\mathbb{R}^{d}}}].
\]
Let $f:x\mapsto e^{\langle\theta,x\rangle_{\mathbb{R}^{d}}}$. This is a
$C^{2}$ function with $\frac{\partial^{2}f}{\partial x_{i}\partial x_{j}%
}=\theta_{i}\theta_{j}f$.

We first need to check the integrability assumption on $f$ in Theorem
\ref{gs-thm}. This is equivalent to $E[e^{\langle\theta,F\rangle
_{\mathbb{R}^{d}}}]<\infty$. To prove this integrability, we compute%
\[
\Gamma_{\langle\theta,F\rangle,\langle\theta,F\rangle}=\sum_{i,j}\theta
_{i}\theta_{j}\Gamma_{ij},
\]
and we note by the positivity of $C-\Gamma$ that this is bounded above almost
surely by the non-random positive constant $K:=\sum_{i,j}\theta_{i}\theta
_{j}C_{ij}$. This implies (see for instance \cite{V}) that $P\left[
\left\langle \theta,F\right\rangle /K>x\right]  \leqslant\Phi\left(  x\right)
$ where $\Phi$ is the standard normal tail. The finiteness of $E[e^{\langle
\theta,F\rangle_{\mathbb{R}^{d}}}]$ follows immediately.

Next, by the positivity of $C-\Gamma$,%
\[
\sum_{i,j}\frac{\partial^{2}f}{\partial x_{i}\partial x_{j}}\left(  x\right)
\left(  \Gamma_{ij}-C_{ij}\right)  =f\left(  x\right)  \sum_{i,j}\theta
_{i}\theta_{j}\left(  \Gamma_{ij}-C_{ij}\right)  \leqslant0.
\]
This is condition (\ref{gs-condition}), so that Theorem \ref{gs-thm} implies
that $E[e^{\langle\theta,F\rangle_{\mathbb{R}^{d}}}]\leqslant E[e^{\langle
\theta,G\rangle_{\mathbb{R}^{d}}}]$ with $G$ a centered Gaussian vector with
covariance matrix $C$. Therefore, since $E[e^{\langle\theta,G\rangle
_{\mathbb{R}^{d}}}]=e^{\frac{1}{2}\langle\theta,C\theta\rangle_{\mathbb{R}%
^{d}}}$, we have
\[
P[F_{1}\geq x_{1},\ldots,F_{d}\geq x_{d}]\leqslant e^{-\langle\theta
,x\rangle_{\mathbb{R}^{d}}+\frac{1}{2}\langle\theta,C\theta\rangle
_{\mathbb{R}^{d}}}\leqslant e^{-\langle\theta,x\rangle_{\mathbb{R}^{d}}%
\frac{1}{2}\Vert C\Vert_{op}\,\Vert\theta\Vert_{\mathbb{R}^{d}}^{2}}.
\]
The desired conclusion follows by choosing $\theta=x/\Vert C\Vert_{op}$, which
represents the optimal choice.\nopagebreak\hspace*{\fill} {\vrule width6pt
height6ptdepth0pt}\bigskip

\section{Universality of the Sherrington-Kirkpatrick model with correlated
media\label{SK}}

Let $N$ be a positive integer, and let $S_{N}=\left\{  -1,1\right\}  ^{N}$,
which represents the set of all possible configurations of the spins of
particles sitting at the integer positions from $1$ to $N$. A parameter
$\beta>0$ is interpreted as the system's inverse temperature. Denote by
$d\sigma$ the uniform probability measure on $S_{N}$, i.e. such that for every
$\sigma\in S_{N}$, the mass of $\left\{  \sigma\right\}  $ is $2^{-N}$. For
any Hamiltonian $H$ defined on $S_{N}$, we can define a probability measure
$P_{N}^{H}$ via $P_{N}^{H}\left(  d\sigma\right)  =d\sigma\exp\left(  -\beta
H\left(  \sigma\right)  \right)  /Z_{N}^{H}$ where $Z_{N}^{H}$ is a
normalizing constant. Therefore,
\begin{equation}
Z_{N}^{H}=2^{-N}\sum_{\sigma\in S_{N}}\exp\left(  -\beta H\left(
\sigma\right)  \right)  .\label{Zee}%
\end{equation}
The measure $P_{N}^{H}$ is the distribution of the system's spins under the
influence of the Hamiltonian $H$. The classical Sherrington-Kirkpatrick (SK,
for short) model for spin systems is a random probability measure in which the
Hamiltonian is random, because of the presence of an external random field
$J=\left\{  J_{i,j}:i,j=1,\cdots,N;i>j\right\}  $ where the random variables
$J_{i,j}$ are IID standard normal (and for notational convenience we assume
the matrix $J$ is defined as being symmetric), and $H=H_{N}$ is given by%
\begin{equation}
H_{N}\left(  \sigma\right)  :=\frac{1}{\sqrt{2N}}\sum_{i\neq j}\sigma
_{i}\sigma_{j}J_{i,j}.\label{HN}%
\end{equation}
The fact that the $J_{i,j}$'s are IID implies that there is no geometry in the
spin system. Indeed, in the sense of distributions w.r.t. the law of $J$, the
interactions between the sites $\left\{  1,\cdots,N\right\}  $ implied by the
definition of $P_{N}^{H}$ do not distinguish between how far apart the sites
are. Such a model is usually called \textquotedblleft
mean-field\textquotedblright, for this lack of geometry. The centered Gaussian
character of the external field $J$ is also an important element in the SK
model's definition, particularly because it implies a behavior for $H_{N}$ of
order $\sqrt{N}$, which can be observed for instance by computing the variance
of $H_{N}\left(  \sigma\right)  $ w.r.t. $J$ for any fixed spin configuration
$\sigma$: it equals $N-1$. A quantity of importance in the study of the
behavior of the measure $P_{N}^{H}$ is its partition function, or free energy,
the scalar $Z_{N}^{H}$ in (\ref{Zee}). In particular, one would like to prove
that it has an almost-sure Lyapunov exponent, namely, a.s. the following limit
exists and is finite:%
\begin{equation}
p\left(  \beta\right)  :=\lim_{N\rightarrow\infty}\frac{1}{N}\log Z_{N}%
^{H}.\label{LyapSK}%
\end{equation}
A proof strategy was defined by Guerra and Toninelli \cite{GT}. In this
classical case, the limit, which we denote by $p_{SK}\left(  \beta\right)  $,
is also known as the Parisi formula (see \cite{MPV} and \cite[page
251]{bovier}). A universality result, where the Gaussian assumption can be
dropped in favor of requiring only three moments for $J$, with the same Parisi
formula for the limit of the normalized log free energy, was established in
\cite{CH}.

In the theorem below, we show that the existence and finiteness of $p\left(
\beta\right)  $, and its equality with $p_{SK}\left(  \beta\right)  $, extends
to external fields $J$ on Wiener space which contain some non-stationarity and
some dependence. Our proof's idea is to use the same smart-path techniques on
Wiener space used in the proofs of Theorems \ref{fs-thm} and \ref{gs-thm}, and
compare $Z_{N}^{H}$ with the free energy of a spin system with IID media
$J^{\ast}$. As explained in more detail in Remark \ref{SKunivRem} below,
Condition (ii) in the theorem is designed to allow for correlations in $J$,
while Condition (iii) implies that the two random media have some asymptotic
proximity in law.

\begin{thm}
\label{SKunivThm}Let $J=\left\{  J_{i,j}:1\leqslant j<i\right\}  $ and
$J^{\ast}=\left\{  J_{i,j}^{\ast}:1\leqslant j<i\right\}  $ be two families of
centered r.v.'s in\textbf{ }$\mathbb{D}^{1,2}$ such that

\begin{description}
\item[(i)] $\left\{  J_{i,j}^{\ast}:1\leqslant j<i\right\}  $ are IID with
variance 1 and $\Gamma_{J_{i,j}^{\ast},J_{k,\ell}^{\ast}}=0$ for all
$(i,j)\neq(k,\ell)$,

\item[(ii)] $\sum_{1\leqslant j<i\leqslant N}\mathbf{E}\left[  \left\vert
\Gamma_{J_{i,j},J_{k,\ell}}\right\vert \right]  =o\left(  N^{2}\right)  ,$

\item[(iii)] $\sum_{1\leqslant j<i\leqslant N}\mathbf{E}\left[  \left\vert
\Gamma_{J_{i,j},J_{i,j}}-\Gamma_{J_{i,j}^{\ast},J_{i,j}^{\ast}}\right\vert
\right]  =o\left(  N^{2}\right)  ,$

\item[(iv)] $\Gamma_{J_{i,j},J_{k,\ell}^{\ast}}=\Gamma_{J_{i,j}^{\ast
},J_{k,\ell}}=0$ for all $i,j,k,\ell$.
\end{description}

Let $Z_{N}^{H}$ be the free energy relative to $J$, as in (\ref{Zee}),
(\ref{HN}). We have $\lim_{N\rightarrow\infty}N^{-1}\log Z_{N}^{H}%
=p_{SK}\left(  \beta\right)  $ in probability. If moreover there exists
$\varepsilon>0$ such that

\begin{description}
\item[(v)] $\sup_{i,j}\mathbf{E}\left[  \left\vert \Gamma_{J_{i,j};J_{i,j}%
}\right\vert ^{1+\varepsilon}\right]  =:M<\infty$,
\end{description}

then the convergence holds almost surely; more specifically, for any
$\delta<2^{-1}\varepsilon/\left(  1+\varepsilon\right)  $, as $N\rightarrow
\infty$, a.s.%
\[
\frac{1}{N}\log Z_{N}^{H}=p_{SK}\left(  \beta\right)  +o(N^{-\delta}).
\]

\end{thm}

\begin{remarks}
\label{SKunivRem} \textrm{{ }}

\begin{enumerate}
\item \textrm{The model in the theorem is the classical SK model (where $J$ is
IID standard normal) as soon as $\Gamma_{J_{i,j},J_{i,j}}\equiv1$ a.s.}

\item \textrm{The classical universality result of Carmona and Hu in \cite{CH}
assumes that $J$ is IID and has three moments. Here we do away with the IID
assumption for $J$, comparing it to an IID $J^{\ast}$ with two moments,
obtaining new SK-universality classes. }

\item \textrm{Condition (ii) above is a way to control the correlations of
$J$. For instance, it is satisfied as soon as $\mathbf{E}\left[  \left\vert
\Gamma_{J_{i,j},J_{k,\ell}}\right\vert \right]  \leqslant\left(  \left\vert
i-k\right\vert +\left\vert j-\ell\right\vert \right)  ^{-r}$ for $r>2$. Since
by formula (\ref{parts}), $\mathbf{E}\left[  \left\vert \Gamma_{J_{i,j}%
,J_{k,\ell}}\right\vert \right]  \geq\left\vert \mathbf{E}\left[
\Gamma_{J_{i,j},J_{k,\ell}}\right]  \right\vert =\left\vert \mathbf{E}\left[
J_{i,j}J_{k,\ell}\right]  \right\vert $, this implies a corresponding
decorrelation rate. }

\item \textrm{Condition (iii) in this corollary can be understood as a kind of
Cesaro-type convergence in distribution. For illustrative purposes, consider
the case where the comparison is with the SK model: we have $\Gamma
_{J_{i,j}^{\ast},J_{i,j}^{\ast}}\equiv1$, and the interpretation of Condition
(iii) can be made more precise. Indeed, by Theorem 5.3.1 in \cite{NPbook},
this type of convergence roughly leads to convergence of $J_{i,j}$ to a
standard normal as $i$ and/or $j\rightarrow\infty$ with $N$. }
\end{enumerate}

\textrm{ }
\end{remarks}

\emph{Proof of Theorem \ref{SKunivThm}}:\vspace*{0.15in}

\noindent\emph{Step 1: a generic result}. We begin by showing a precursor
result for convergence in probability, for a generic situation. Assume that
$J$ and $J^{\ast}$ satisfy merely (ii), (iii), and (iv). We will show that for
any $f\in C^{2}\left(  \mathbb{R}\right)  $ with $\Vert f^{\prime}%
\Vert_{\infty}\leqslant1$ and $\Vert f^{\prime\prime}\Vert_{\infty}\leqslant
1$,
\begin{equation}
\left\vert \mathbf{E}\left[  f(\frac{1}{N}\log Z_{N}^{\ast H})\right]
-\mathbf{E}\left[  f(\frac{1}{N}\log Z_{N}^{H})\right]  \right\vert =o(1).
\label{generic}%
\end{equation}

We compactify the notation by reindexing the set $\left\{
i,j:i>j;i,j=1,\cdots,N\right\}  $ as the set $\left\{  1,2,\cdots,\bar
{N}\right\}  $ where $\bar{N}:=N(N-1)/2$, with a bijection mapping each
$n=1,\cdots,\bar{N}$ to a pair $\left(  i,j\right)  $, using any fixed
bijection, with $\bar{J}_{n}:=J_{i,j}$, $\bar{J}_{n}^{\ast}:=J_{i,j}^{\ast}$,
and $\tau_{n}:=\sigma_{i}\sigma_{j}$, with $P_{\sigma}$ the uniform
probability measure on $S_{N}$, so that each r.v. $\tau_{n}$ under $P_{\sigma
}$ is dominated by $1$. We use $\bar{J}$ and $\bar{J}^{\ast}$ to denote the
corresponding $\bar{N}$-dimensional random vectors.

Fix $\gamma>0$, $c\in\lbrack0,1]$ and $f$ as above. We define for any vector
$u\in\mathbb{R}^{\bar{N}}$, and $t\in\lbrack0,1]$,
\begin{align*}
Z(\gamma,u)  &  :=E_{\sigma}\left[  \exp\left(  \gamma\sum_{n=1}^{\bar{N}}%
\tau_{n}u_{n}\right)  \right]  ,\\
\varphi(t)  &  :=\mathbf{E}[f(c\log Z(\gamma,\sqrt{t}\bar{J}^{\ast}+\sqrt
{1-t}\bar{J}))].
\end{align*}
For $i=1,\ldots,\bar{N}$ and $u\in\mathbb{R}^{\bar{N}}$, we define%
\[
h_{i}(u):=\frac{E_{\sigma}[\tau_{i}e^{\gamma\sum_{n=1}^{\bar{N}}\tau_{n}u_{n}%
}]}{E_{\sigma}[e^{\gamma\sum_{n=1}^{\bar{N}}\tau_{n}u_{n}}]}\,f^{\prime}(c\log
E_{\sigma}[e^{\gamma\sum_{n=1}^{\bar{N}}\tau_{n}u_{n}}]).
\]
We compute that for any $i,j=1,\ldots,\bar{N}$, we have $\frac{\partial h_{i}%
}{\partial u_{j}}(u)=\gamma\,S_{i,j}(u)$ where%
\begin{align*}
S_{i,j}(u)  &  :=\left(  \frac{E_{\sigma}[\tau_{i}\tau_{j}e^{\gamma\sum
_{n=1}^{\bar{N}}\tau_{n}u_{n}}]}{E_{\sigma}[e^{\gamma\sum_{n=1}^{\bar{N}}%
\tau_{n}u_{n}}]}-\frac{E_{\sigma}[\tau_{i}e^{\gamma\sum_{n=1}^{\bar{N}}%
\tau_{n}u_{n}}]E_{\sigma}[\tau_{j}e^{\gamma\sum_{n=1}^{\bar{N}}\tau_{n}u_{n}%
}]}{E_{\sigma}[e^{\gamma\sum_{n=1}^{\bar{N}}\tau_{n}u_{n}}]^{2}}\right)
\,f^{\prime}(c\log E_{\sigma}[e^{\gamma\sum_{n=1}^{\bar{N}}\tau_{n}u_{n}}])\\
&  +c\,\frac{E_{\sigma}[\tau_{i}e^{\gamma\sum_{n=1}^{\bar{N}}\tau_{n}u_{n}%
}]E_{\sigma}[\tau_{j}e^{\gamma\sum_{n=1}^{\bar{N}}\tau_{n}u_{n}}]}{E_{\sigma
}[e^{\gamma\sum_{n=1}^{\bar{N}}\tau_{n}u_{n}}]^{2}}\,f^{\prime\prime}%
(c\times\log E_{\sigma}[e^{\gamma\sum_{n=1}^{\bar{N}}\tau_{n}u_{n}}]).
\end{align*}
Notice that since $c,\tau_{i},f^{\prime}$, and $f^{\prime\prime}$ are all
dominated by $1$, we get $\left\vert S_{i,j}(u)\right\vert \leqslant3$. Using
the chain rule of standard calculus,%
\[
\varphi^{\prime}(t)=\frac{c\,\gamma}{2}\,\sum_{i=1}^{\bar{N}}\left\{  \frac
{1}{\sqrt{t}}\mathbf{E}[\bar{J}_{i}^{\ast}h_{i}(\sqrt{t}\bar{J}^{\ast}%
+\sqrt{1-t}\bar{J})]-\frac{1}{\sqrt{1-t}}\mathbf{E}[\bar{J}_{i}h_{i}(\sqrt
{t}\bar{J}^{\ast}+\sqrt{1-t}\bar{J})]\right\}  .
\]
Now using the integration-by-parts formula on Wiener space (\ref{parts}), and
Condition (iv), this computes as%
\begin{align*}
\varphi^{\prime}(t)  &  =\frac{c\,\gamma^{2}}{2}\sum_{i=1}^{\bar{N}}%
\mathbf{E}\left[  S_{i,i}(\sqrt{t}\bar{J}^{\ast}+\sqrt{1-t}\bar{J}%
)(\Gamma_{\bar{J}_{i}^{\ast},\bar{J}_{i}^{\ast}}-\Gamma_{\bar{J}_{i},\bar
{J}_{i}})\right] \\
&  +\frac{c\,\gamma^{2}}{2}\sum_{1\leqslant i\neq j\leqslant\bar{N}}%
\mathbf{E}\left[  S_{i,j}(\sqrt{t}\bar{J}^{\ast}+\sqrt{1-t}\bar{J}%
)\Gamma_{\bar{J}_{i},\bar{J}_{j}}\right] .
\end{align*}
The boundedness of $\left\vert S_{i,j}(u)\right\vert $ by $3$ yields, by
integrating over $t\in\lbrack0,1]$, that%
\begin{align*}
&  \left\vert \mathbf{E}\left[  f(c\log Z_{\bar{N}}(\gamma,\bar{J}^{\ast
}))\right]  -\mathbf{E}\left[  f(c\log Z_{\bar{N}}(\gamma,\bar{J}))\right]
\right\vert =\left\vert \int_{0}^{1}\varphi^{\prime}(t)dt\right\vert \\
&  \leqslant\frac{3c\,\gamma^{2}}{2}\sum_{i=1}^{\bar{N}}\mathbf{E}\left[
\left\vert \Gamma_{\bar{J}_{i}^{\ast},\bar{J}_{i}^{\ast}}-\Gamma_{\bar{J}%
_{i},\bar{J}_{i}}\right\vert \right]  +\frac{3c\,\gamma^{2}}{2}\sum
_{1\leqslant i\neq j\leqslant\bar{N}}\mathbf{E}\left[  \left\vert \Gamma
_{\bar{J}_{i},\bar{J}_{j}}\right\vert \right]  .
\end{align*}
By Conditions (ii) and (iii), replacing $\gamma$ by $\beta/\sqrt{N}$ and $c$
by $1/N$, with $\bar{N}=N(N-1)/2$, relation (\ref{generic}) follows.

\bigskip

\noindent\emph{Step 2: Convergences.} In this step we assume for the moment
that $\lim_{N\rightarrow\infty}N^{-1}\log Z_{N}^{\ast H}=p_{SK}\left(
\beta\right)  $ holds in probability. This convergence is established below in
Step 3. Combining this convergence and relation (\ref{generic}), we get that
$N^{-1}\log Z_{N}^{H}$ converges in distribution, and thus in probability, to
$p_{SK}\left(  \beta\right)  $, which is the first conclusion of the theorem.
To establish the second conclusion, i.e. the almost-sure convergence, let%
\[
F_{N}:=\frac{1}{N}\log Z_{N}^{H} - \frac{1}{N}\mathbf{E}\left[  \log Z_{N}%
^{H}\right] .
\]
By the chain rule of Malliavin calculus, and using the notation $E_{N}^{H}$
for expectations of functions of the configuration $\sigma$ under the polymer
measure defined by
\[
P_{N}^{H}\left(  \left\{  \sigma\right\}  \right)  =\frac{1}{2^{N}}\frac
{\exp\left(  -\beta~H_{N}\left(  \sigma\right)  \right)  }{\sum_{\sigma\in
S_{N}}\exp\left(  -\beta~H_{N}\left(  \sigma\right)  \right)  },
\]
we compute
\begin{align*}
DF_{N}  &  =\frac{1}{N}\frac{1}{Z_{N}^{H}}\left(  -\beta2^{-N}\right)
\sum_{\sigma\in S_{N}}\exp\left(  -\beta~H_{N}\left(  \sigma\right)  \right)
DH_{N}\left(  \sigma\right) \\
&  =\frac{-\beta}{N}E_{N}^{H}\left[  DH_{N}\left(  \sigma\right)  \right]  .
\end{align*}
Now, using the intermediary of the Mehler formula (see, e.g.,
\cite[Proposition 3.7]{nourdin-viens}), it is easy to check that we can
express%
\[
\Gamma_{F_{N},F_{N}}=\frac{\beta^{2}}{N^{2}}E_{N}^{H}\otimes\tilde{E}_{N}%
^{H}\left[  \Gamma_{H_{N}\left(  \sigma\right)  ,H_{N}\left(  \tilde{\sigma
}\right)  }\right]
\]
where for fixed random medium $J$, under $P_{N}^{H}\otimes\tilde{P}_{N}^{H}$,
$\left(  \sigma,\tilde{\sigma}\right)  $ are two independent copies of
$\sigma$ under the polymer measure $P_{N}^{H}$. We compute for any
$\sigma,\sigma^{\prime}\in S_{N}$,%
\[
\Gamma_{H_{N}\left(  \sigma\right)  ,H_{N}\left(  \sigma^{\prime}\right)
}=\frac{2}{N}\sum_{1\leqslant j<i\leqslant N}\Gamma_{J_{i,j},J_{i,j}}%
\sigma_{i}\sigma_{i}^{\prime}\sigma_{j}\sigma_{j}^{\prime}.
\]
Since $\left\vert \sigma_{i}\right\vert =1$ for any $\sigma\in S_{N}$, we get
\begin{equation}
\left| \Gamma_{F_{N},F_{N}}\right| \leqslant\frac{2\beta^{2}}{N^{3}}%
\sum_{1\leqslant j<i\leqslant N}\left\vert \Gamma_{J_{i,j},J_{i,j}}\right\vert
. \label{GammaFbound}%
\end{equation}
By Assumption (v), $\mathbf{E}\left[  \left\vert \Gamma_{J_{i,j},J_{i,j}%
}\right\vert ^{1+\varepsilon}\right]  $ is uniformly bounded by $M$.
Therefore, using Jensen's inequality for the uniform measure on the set
$\left\{  i,j=1,\cdots,N;~i>j\right\}  $ and the power function $\left\vert
x\right\vert ^{1+\varepsilon}$,
\begin{align*}
\mathbf{E}\left[  \left\vert \Gamma_{F_{N},F_{N}}\right\vert ^{1+\varepsilon
}\right]   &  \leqslant\left(  \frac{2\beta^{2}}{N^{3}}\right)
^{1+\varepsilon}\left(  N(N-1)/2\right)  ^{1+\varepsilon}\frac{2}{N(N-1)}%
\sum_{1\leqslant j<i\leqslant N}\mathbf{E}\left[  \left\vert \Gamma
_{J_{i,j},J_{i,j}}\right\vert ^{1+\varepsilon}\right] \\
&  \leqslant M\beta^{2+2\varepsilon}N^{-1-\varepsilon}.
\end{align*}
We now need a Poincar\'{e}-type inequality on Wiener space relative to the
operator $\Gamma$, which is recorded and proved below in Lemma
\ref{PoincareLem}: applying this lemma with $F=F_{N}$ and $p=2+2\varepsilon$
yields%
\[
\mathbf{E}\left[  \left\vert F_{N}\right\vert ^{2+2\varepsilon}\right]
\leqslant\left(  1+2\varepsilon\right)  ^{1+\varepsilon}M\beta^{2+2\varepsilon
}N^{-1-\varepsilon}.
\]
A standard application of the Borel-Cantelli lemma via Chebyshev's inequality
yields that for any $\delta<2^{-1}\varepsilon/\left(  1+\varepsilon\right)  $,
almost surely, $F_{N}=o(N^{-\delta})$, as announced in the theorem.\vspace
*{0.1in}

\noindent\emph{Step 3: Conclusion}. To finish the proof of the theorem, we
only need to show that $\lim_{N\rightarrow\infty}N^{-1}\log Z_{N}^{\ast
H}=p_{SK}\left(  \beta\right)  $ holds in probability. The universality result
of Carmona and Hu as stated in \cite{CH} shows that this convergence holds if
we assumed in addition that $J_{i,j}^{\ast}$ had a finite third moment.
However, an inspection of their proof reveals that the convergence holds in
probability without the third moment condition: one may use a computation
similar to the calculation in Step 1 above, to establish this; the details are
omitted. \nopagebreak\hspace*{\fill} {\vrule width6pt height6ptdepth0pt}%
\bigskip

\begin{lemme}
\label{PoincareLem}For any centered $F\in\mathbf{D}^{1,2}$, and any $p\geq2$,%
\[
\mathbf{E}\left[  \left\vert F\right\vert ^{p}\right]  \leqslant\left(
p-1\right)  ^{p/2}\mathbf{E}\left[  \left\vert \Gamma_{F,F}\right\vert
^{p/2}\right]  .
\]

\end{lemme}

\noindent\textit{Proof}. For $p=2$, by relation (\ref{parts}), the inequality
holds almost as an equality (one has $\mathbf{E}[F^{2}]=\mathbf{E}%
[\Gamma_{F,F}]\leqslant\mathbf{E}[|\Gamma_{F,F}|]$). Therefore we assume
$p>2$. With the notation $G\left(  x\right)  =\mathrm{sgn}\left(  x\right)
\left\vert x\right\vert ^{p-1}$, and thus $G^{\prime}\left(  x\right)
=\left(  p-1\right)  \mathrm{sgn}\left(  x\right)  \left\vert x\right\vert
^{p-2}$, and $G\left(  F_{N}\right)  \in\mathbf{D}^{1,2}$ with $D\left(
G\left(  F\right)  \right)  =\left(  p-1\right)  \mathrm{sgn}\left(  F\right)
\left\vert F\right\vert ^{p-2}DF$, we have, using again (\ref{parts}),
\[
\mathbf{E}\left[  \left\vert F\right\vert ^{p}\right]  =\mathbf{E}\left[
FG\left(  F\right)  \right]  =\left(  p-1\right)  \mathbf{E}\left[
\mathrm{sgn}\left(  F\right)  \left\vert F\right\vert ^{p-2}\Gamma
_{F,F}\right]  .
\]
Now invoking H\"{o}lder's inequality we get%
\[
\mathbf{E}\left[  \left\vert F\right\vert ^{p}\right]  \leqslant\left(
p-1\right)  \mathbf{E}\left[  \left\vert \Gamma_{F,F}\right\vert
^{p/2}\right]  ^{2/p}\mathbf{E}\left[  \left\vert F\right\vert ^{p}\right]
^{1-2/p}.
\]
The lemma follows immediately.\nopagebreak\hspace*{\fill} {\vrule width6pt
height6ptdepth0pt}


\bigskip

\begin{thebibliography}{99}                                                                                               %


\bibitem {adler}R. J. Adler (1990). \emph{An introduction to continuity,
extrema, and related topics for general Gaussian processes} Lecture
Notes--Monograph Series \textbf{12}, Hayward, CA. Institute of Mathematical Statistics.

\bibitem{adlertaylor} R.$\,$J. Adler and J.E. Taylor (2007). {\it Random fields and geometry. }
Springer-Verlag.

\bibitem {AMV}H. Airault, P. Malliavin, F. Viens (2010). Stokes formula on the
Wiener space and $n$-dimensional Nourdin-Peccati analysis. \emph{J. Funct.
Anal.} \textbf{258} (5), 1763-1783

\bibitem {anderson}T.W. Anderson (1955). The integral of a symmetric unimodal
function over a symmetric convex set and some probability inequalities.
\textit{Proc. Amer. Math. Soc.} \textbf{6}, 170-176.

\bibitem {bovier}A. Bovier (2006). \emph{Statistical mechanics of disordered
systems. A mathematical perspective.} Cambridge University Press.

\bibitem {CH}Ph. Carmona, Y. Hu (2006). Universality in
Sherrington-Kirkpatrick's spin glass model. \emph{Annales IHP (B) Prob. Stat.}
\textbf{42} (2), 215-222.

\bibitem{Chatterjee}S. Chatterjee (2005). An error bound in the Sudakov-Fernique inequality. ArXiv:math/0510424.

\bibitem {GT}F. Guerra, F.L. Toninelli (2002). The thermodynamic limit in mean
eld spin glass models. \emph{Comm. Math. Phys}. \textbf{230} , no. 1, 71-79.

\bibitem {MPV}M. M\'{e}zard, G. Parisi, and M. A. Virasoro (1987). \emph{Spin
Glass Theory and Beyond}, World Scientific Lecture Notes in Physics, vol.
\textbf{9}. World Scientific.

\bibitem {nourdin-peccati}I. Nourdin and G. Peccati (2009). Stein's method on
Wiener chaos. \textit{Probab. Theory Related Fields} \textbf{145}, 75-118.

\bibitem {NPbook}I. Nourdin, G. Peccati (2012). \emph{Normal approximation
with Malliavin calculus: from Stein's method to universality}. Cambridge
University Press.

\bibitem {NPR}I. Nourdin, G. Peccati, A. R\'{e}veillac (2010). Multivariate
normal approximation using Stein's method and Malliavin calculus. \emph{Ann.
IHP (B) Probab. Statist.} \textbf{46} (1), 45-58.

\bibitem {nourdin-simon}I. Nourdin and T. Simon (2006). On the absolute
continuity of one-dimensional SDEs driven by a fractional Brownian motion.
\textit{Stat. Probab. Lett.} \textbf{76}, no. 9, 907-912.

\bibitem {nourdin-viens}I. Nourdin and F.G. Viens (2009). Density formula and
concentration inequalities with Malliavin calculus. \textit{Electron. J.
Probab.} \textbf{14}, 2287-2309.

\bibitem {Nbook}D. Nualart (2006). \emph{Malliavin calculus and related
topics.} Springer Verlag.

\bibitem {nualart-ouknine}D. Nualart and Y. Ouknine (2002). Regularization of
differential equations by fractional noise. \textit{Stoch. Proc. Appl.}
\textbf{102}, no. 1, 103-116.

\bibitem {nualart-saussereau}D. Nualart and B. Saussereau (2009). Malliavin
calculus for stochastic differential equations driven by a fractional Brownian
motion. \textit{Stoch. Proc. Appl.} \textbf{119}, no. 2, 391-409.

\bibitem {piterbarg}V.I. Piterbarg (1982). Gaussian random processes.
[Progress in Science and Technology] \textit{Teor. Veroyatnost. Mat. Statist.
Teor Kibernet.} \textbf{9}, 155-198.

\bibitem {slepian}D. Slepian (1962). The one-sided barier problem for Gausian
noise. \textit{Bell. Syst. Tech. J.} \textbf{41}, no. 2, 463-501.

\bibitem {U}A.-S. \"{U}st\"{u}nel (1995). \emph{An introduction to analysis on
Wiener space}. Springer Verlag.

\bibitem {V}F. Viens (2009). Stein's lemma, Malliavin calculus, and tail
bounds, with application to polymer fluctuation exponent. \emph{Stochastic
Processes and their Applications} \textbf{119}, 3671-3698.

\bibitem {Vitale}R.A. Vitale (2000). Some comparisons for Gaussian processes.
\textit{Proc. Amer. Math. Soc.} \textbf{128}, 3043-3046.
\end{thebibliography}
\end{document}